\documentclass[10pt,twoside]{article}
\usepackage[utf8]{inputenc}
\usepackage{natbib}
\usepackage{graphicx}
\usepackage{subfloat}
\usepackage{subfig}
\usepackage{booktabs}
\usepackage{adjustbox}
\usepackage[table,xcdraw]{xcolor}
\usepackage{amsmath}
\DeclareMathOperator*{\argmax}{\arg\!\max}

\DeclareMathOperator*{\med}{med}
\usepackage{amssymb}
\usepackage[]{todonotes}
\usepackage{fancyhdr}
\usepackage{multirow}
\usepackage{hyperref}

\usepackage[english]{babel}

\usepackage[none]{hyphenat}

\usepackage[Alg.]{algorithm}
\usepackage{algpseudocode}

\newcommand{\ignore}[1]{}

\title{A robust approach for principal component analysis}
\author{\emph{María Camila Vásquez Correa} and \emph{Henry Laniado Rodas}\\
\vspace{0.3cm}
\small{\tt{mvasqu49@eafit.edu.co, hlaniado@eafit.edu.co}}\\
Mahematical sciences department\\
School of sciences\\
Universidad EAFIT\\
Medellín -- Colombia}
\date{}

\usepackage{anysize}
\marginsize{3cm}{2cm}{2cm}{3cm}

\begin{document}
\maketitle
\begin{abstract}
In this paper we analyze different ways of performing principal component analysis throughout three different
approaches: robust covariance and correlation matrix estimation, projection pursuit approach and
non-parametric maximum entropy algorithm. The objective of these approaches is the correction of the
well known sensitivity to outliers of the classical method for principal component analysis. Due to their
robustness, they perform very well in contaminated data, while the classical approach fails to preserve
the characteristics of the core information.
\vspace{0.2cm}\\
\noindent \textbf{Keywords:}
Statistics, non-parametric, robust, PCA.
\end{abstract}
\section{Introduction}
In principal component analysis (PCA), we seek to maximize the variance of a linear combination of a set of independent variables considered \citep{rencher2003methods}. Essentially, the goal of PCA is to identify the most meaningful basis to re-express a data set in a way that it preserves its structure and leaves behind some of its noise. \\[5pt]
As many other techniques, PCA has been proved to be sensitive to outliers in the data set by various authors \citep{robustPCA}. To face this problem, there are several courses of action. The first, is to apply a robust estimator of the covariance or correlation matrix, to give full weight to observations assumed to come from the main body of the data, but reduced weight to the tails of the contaminated distribution. \\[5pt]
Let $X = (X_1, X_2, \dots, X_p)$ be a multivariate vector with distribution function $F$ in $\mathbb{R}^p$. Let, also, $\Sigma$ be the correlation matrix of $X$, the first eigenvector of this matrix is defined as a unit length vector, $\hat{v}_1$, which maximizes the dispersion of the projection of the observation on that direction. The second eigenvector is defined similarly, but now we only maximize over all vectors perpendicular to the first eigenvector. The k-th eigenvector $\hat{v}_k$ is defined as the maximizer of the function
\begin{equation}
\label{eq:pca}
    a \rightarrow S(a^tx_1, \dots, a^tx_n)
\end{equation}
under the restrictions
\begin{align}
    a\bot \hat{v}_1, \dots a\bot \hat{v}_{k-1}, & & a^ta = 1
\end{align}
Where $S$ is a function of the dispersion of the sample. The corresponding eigenvalues are given by
\begin{equation}
    \hat{\lambda}_k =S(\hat{v}_k^tx_1, \dots, \hat{v}_k ^tx_n)
\end{equation}
for $1 \leq k \leq n$. In classical PCA one takes for $S$ in Equation \ref{eq:pca} the square root of the sample variance, and the solutions to the above problem are given by the eigenvectors and eigenvalues of the sample covariance matrix. This measure, the sample covariance matrix, is very sensitive to outliers. Some approaches to PCA is to calculate this matrix using robust measures.\\[5pt]
In this project we focus on some ways of calculating this matrix, and compare it to some results in the literature like the projection pursuit purposed in \cite{croux1996fast} and the maximum entropy approach purposed in \cite{maxEnt}, all of these approaches to principal component analysis robust to outliers.

\section{Robust PCA}

Consider the solution to the optimization problem in \ref{eq:pca} given 
by the eigen values of the covariance matrix $\Sigma$, or the correlation
matrix for the standarized data. Our first approach to perform robust PCA
is to generate a robust estimator to $\Sigma$, such as the sensitivity to
outliers is properly corrected.

\subsection{Robust covariance matrix by M-estimator}
Let $\bar{X}$ represents the $p\times 1$ vector of sample means and $V$
the covariance matrix, then, the Mahalanobis squared distance of the
$m$th observation form the mean of the observations is defined by
\begin{align}
    \label{mahal}
    d_m^2 & = (x_m - \bar{X})^TV^{-1}(x_m-\bar{X}) & & m = 1, \dots n
\end{align}
Atypical multivariate vectors of observations will tend to deflate
correlations and possibly inflate variances, an this will decrease the
Mahalanobis distance for the outliers.\\[5pt]
Based on this approach, one can propose a M-estimator for the sample
mean and variance, that gives a modification of the classical estimators,
but give full weight to observations assumed to come from the main
body of the data and reduced weight to possible outliers, that is, 
observations with large Mahalanobis distance.\\[5pt]
\cite{robustP} define the robust estimators for mean and covariance as
follows
\begin{align}
    \label{eq:robMean}
    \hat{\mathbf{x}} & = \frac{\sum_{m = 1}^n w_mx_m}{\sum_{m = 1}^n w_m}\\
    \label{eq:robVar}
    \mathbf{V} & = 
    \frac{\sum_{m = 1}^n w_m^2(x_m - \bar{X})(x_m - \bar{X})^T}
    {\sum_{m = 1}^n w_m^2 - 1}
\end{align}
where.
\begin{itemize}
    \item $w_m = w(d_m) = \omega(d_m)/d_m$
    \item $d_m = \{(x_m - \bar{X})^TV^{-1}(x_m-\bar{X})\}^{1/2}$
\end{itemize}
The solution for \ref{eq:robMean} and \ref{eq:robVar} is iterative and
best described in \citep{robustP}. This estimator can be used instead 
of the usual estimator to calculate the eigenvectors that provide the
directions for the principal components.

\subsection{Robust correlation matrix by Kendall's and Spearman's
correlation coefficients}
Kendall’s Tau and Spearman’s rank correlation coefficient assess
statistical associations based on the ranks of the data. Ranking data 
is carried out on the variables that are separately put in order and are
numbered. \\[5pt]
Correlation coefficients take the values between minus one and plus one.
The positive correlation signifies that the ranks of both the variables
are increasing.  On the other hand, the negative correlation signifies
that as the rank of one variable is increased, the rank of the other
variable is decreased.\\[5pt]
The widely used parametric correlation coefficient, known as the Pearson
product–moment correlation coefficient is defined as
\begin{equation}
    \label{eq:corr}
    \rho = \frac{\sum (X_i - \bar{X})(Y_i - \bar{Y})}
    {\left[\sum (X_i - \bar{X})^2\sum(Y_i - \bar{Y})^2\right]^{1/2}}
\end{equation}
As for the Spearman's rank correlation, \ref{eq:corr} can be re-written
as \citep{spearman}
\begin{equation}
\label{eq:spearman}
    \rho_s = \frac{12\sum\mathcal{R}(X_i)\mathcal{R}(Y_i)}{n^3 - n}
    - \frac{3(n+1)}{n-1}
\end{equation}
where $\mathcal{R}(X_i)$ is the rank of the observation $X_i$. \\[5pt]
For the Kendall's tau, consider $(X_i, Y_i)$ and $(X_j, Y_j)$, a pair of
bivariate observations, If $X_i - X_j$ and $Y_i - Y_j$ have the same
sign, the pair is \textit{concordant}, else, is \textit{discordant}.
Let $C$ be the number of concordant pairs and $D$ the number of
discordant ones. The Kendall's S, $S = C-D$ , measures the strength of
the relationship between two variables. Then, the Kendall's tau is 
defined as follows \citep{kendall}
\begin{equation}
    \label{eq:kendall}
    \hat{\tau} =  \frac{2S}{n(n-1)}
\end{equation}
Using these estimators for the correlations, that are both no parametric,
one can compute the correlation matrix and use it to perform a principal
component analysis on a scaled data set, and expect it to perform in a robust
way.

\section{PCA based on projection pursuit}

\subsection{Introduction}
Projection pursuit (PP) methods aim at finding structures in
multivariate data by projecting them on a lower-dimensional subspace, 
often of dimension one, 
selected by maximizing a certain \textit{projection index} 
\citep{algs}. \\[5pt]
PCA is an example of the PP approach.
If we have $n$ observations, each of them column
vectors of dimension $p$, the first principal component is
obtained by finding the unit vector \textbf{\textit{a}} which
maximizes the variance of the data projected on it:
\begin{equation}
\label{eq:pp}
    \mathbf{a}_1 = \argmax\limits_{
    ||\mathbf{a}|| = 1}\ S^2(\mathbf{a}^t\mathbf{x}_1,
    \dots,\mathbf{a}^t\mathbf{x}_n)
\end{equation}
where $S^2$ is the variance. By projecting the data on the
direction $\mathbf{a}_1$ we obtain univariate data, in
accordance to PP idea. Taking the variance as a projection 
index leads, to standard PCA. But, taking a robust measure of
variance can lead to a robust procedure for PCA. \\[5pt]
If we have already computed the $(k-1)$th principal component,
then the direction of the $k$th component, with $1<k\leq p$, is 
defined as the unit vector maximizing the index $S^2$ of the data 
projected
on it, under the condition of being orthogonal to all
previously obtained components:
\begin{equation}
\label{eq:pp1}
    \mathbf{a}_1 = \argmax\limits_{
    ||\mathbf{a}|| = 1,a\bot a_1 \dots a \bot a_{k-1}}\ S^2(\mathbf{a}^t\mathbf{x}_1,
    \dots,\mathbf{a}^t\mathbf{x}_n)
\end{equation}
following this idea, one can only compute a fraction of the components,
implying reduction in time and space when $p$ is large.

\subsection{Algorithm}
Consider a data matrix $\mathbf{X}$ with $n$ rows and $p$
columns, having the observations $\mathbf{x}_i = (x_{i1}, \dots, x_{in})^t$ in its columns. 
Assume, also, that $p \leq n$.\\[5pt]
The maximal values for the variances, denoted by
\begin{equation}
    \label{eq:varpp}
    \lambda_k = S^2(\mathbf{a}_k^t\mathbf{x}_1,
    \dots,\mathbf{a}_k^t\mathbf{x}_n)
\end{equation}
will be provided by the algorithm and also the directions of the
vectors in \ref{eq:pp} and \ref{eq:pp1}. For the method, we will
consider robust versions of $S^2$: the Median Absolute Deviation
(MAD):
\begin{equation}
    \label{eq:MAD}
    \text{MAD}(z_1,\dots,z_n) = 
    1.48\med\limits_{j}|z_j - \med\limits_{i}z_i|
\end{equation}
and the first quartile of the pairwise differences between all data
points
\begin{equation}
    \label{eq:qn}
    Q(z_1,\dots,z_n) = 
    2.22\{|z_i - z_j|; 
    1\leq i < j \leq n\}_{\left(
    \begin{array}{cc}
        n \\ 2
    \end{array}\right)/4}
\end{equation}
Let $\mathbf{x}_i^1 = \mathbf{x}_i - \hat{\mathbf{\mu}}(\mathbf{X})$
be the centered data, where $\hat{\mathbf{\mu}}(\mathbf{X}$ is
the $L_1$-median, a robust estimator of the center of the data.\\[5pt]
Suppose that in step $k-1$, the algorithm returned the vector
$\hat{\mathbf{a}}_{k-1}$, an approximation for the solution to
\ref{eq:pp1}. Then we update the observations according to
\begin{equation}
    \label{eq:update}
    \mathbf{x}_i^k = \mathbf{x}_i^{k-1} - 
    (\hat{\mathbf{a}}_{k-1}^t\mathbf{x}_i^{k-1})
    \hat{\mathbf{a}}_{k-1}
\end{equation}
for $1\leq i \leq n$ and $k > 1$. The algorithm only considers $n$
trial directions in the set
\begin{equation}
    \label{eq:set}
    A_{n,k} = \left\lbrace
    \frac{\mathbf{x}_1^{k}}{||\mathbf{x}_1^{k}||}, \dots,
     \frac{\mathbf{x}_n^{k}}{||\mathbf{x}_n^{k}||}
    \right\rbrace
\end{equation}
And then, the $k$th eigenvalue is approximated by
\begin{equation}
    \label{eq:eigen}
    \hat{\lambda}_k = \max\limits_{\mathbf{a}\in A_{n,k}} 
    S^2(\mathbf{a}_k^t\mathbf{x}_1,
    \dots,\mathbf{a}_k^t\mathbf{x}_n)
\end{equation}
and $\hat{\mathbf{a}}_{k}$ is the argument where the maximum of
\ref{eq:eigen} is obtained \citep{croux2005high}. This will be, then, 
the direction of the $k$th robust principal component.

\section{PCA based on non-parametric maximum entropy}
Consider a data set of samples $X = (X_1, \dots, X_n)$, where $X_i$ us a variable with dimensionallity $d$, $U = (U_1, \dots, U_m) \in  \mathbb{R}^{d\times m}$ a projection matrix whose columns constitute the bases of a $m$-dimensional subspace, and $V = (V_1, \dots, V_n) \in  \mathbb{R}^{m\times n}$ is the projection coordinates under the projeciton matrix $U$. \\[5pt]
PCA can be formulated as the following optimization problem:
\begin{equation}
\label{eq:2}
    \min\limits_{U,V}\sum_{i = 1}^n ||X_i - (\mu + UV_i)||^2
\end{equation}
where $\mu$ is the center of $X$. The global minimum of~\ref{eq:2} is provided, as we know, by singular value decomposition, whose optimal solution is also the solution of the equation:
\begin{equation}
\label{eq:or}
    \max\limits_{U^TU = I} Tr(U^T\Sigma U)
\end{equation}
where $\Sigma$ is the covariance matrix, $Tr(\dot)$ is the matrix trace operation. Equation~\ref{eq:or} searches for a projection matrix where the variances of $U^TX$ are maximized. \\[5pt]
The aim of MaxEnt-PCA is to learn a new data distribution in a subspace such that entropy is maximized. The Reinyi's quadratic entropy of a random variable $X$ with probability density function $f_X(x)$ defined by
\begin{equation}
    \label{eq:ent}
    H(X) = -\log\int f^2_X(x) dx
\end{equation}

When we replace the density function by its estimator, a Gaussian Kernel
\begin{equation}
    \label{eq:gk}
    G(X - X_i, \sigma^2) =  \frac{1}{(2\pi)^{d/2}}\exp{\left(-\frac{||x-x_i||^2}{2\sigma^2}\right)}
\end{equation}
we obtain the following optimization problem
\begin{align}
\label{eq:maxEnt}
   &  \max\limits_U \left(-\log\left(\frac{1}{n^2}\sum_{i = 1}^n \sum_{j=1}^n G(U^TX_j - U^tX_i, \sigma^2\right)\right) & & \text{s.t.}\quad U^tU = I
\end{align}
where $G(X - X_i, \sigma^2)$ is the Gaussian kernel with bandwidth $\Sigma = \sigma^2 I$\\[5pt]
The optimal solution of MaxEnt-PCA in~\ref{eq:maxEnt} is given by the eigenvectors of the following generalized eigen-decomposition problem:
\begin{equation}
    \label{eq:edp}
    XL(U)X^TU = 2U\Lambda
\end{equation}
where:
\begin{align}
\label{eq:laque}
    L(U) & = D(U) - W(U) \\
    W_{ij}(U) & = \frac{G(U^Tx_i - U^Tx_j, \sigma^2)}{\sigma^2 \sum_{i = 1}^nG(U^Tx_i - U^Tx_j, \sigma^2} \\
    D_{ii} & = \sum_{j = 1}^n W_{ij}U
\end{align}
Since $L(U)$ in \ref{eq:laque} is  also  a  function  of $U$,  the  eigenvalue decomposition problem in  \ref{eq:laque} has  no  closed-form  solution. We use, then a gradient-based fixed-point algorithm, and follow these steps to update the projection matrix $U$:
\begin{align}
    U & = (I + \beta XL(U)X^T)U \\
    U & = svd(U)
\end{align}
where $\beta$ is a step length to ensure an increment of the objective function, and $svd(U)$ returns an orthonormal base by the Singular Value Decomposition (SVD)  on  matrix $U$. Also, the bandwidth $\sigma^2$, as a factor of average distance between points is given by
\citep{maxEnt}
    \begin{equation}
        \sigma^2 = \frac{1}{sn^2}
        \sum_{i=1}^n\sum_{j=1}^n||U^Tx_i - U^Tx_j||^2
    \end{equation}
where $s$ is a scale factor. in each iteration, we want to see how much entropy we have achieved it, and if the difference between iterations has converged. More details about the performance of this method are in \citep{maxEnt}.
\section{Numerical experiments}
\subsection{Robust estimators for the correlation and covariance matrix}
\begin{figure}[ht!]
\centering
    \subfloat[PCA using the Kendall's $ \tau$ for correlation matrix]{\includegraphics[scale = 0.35]{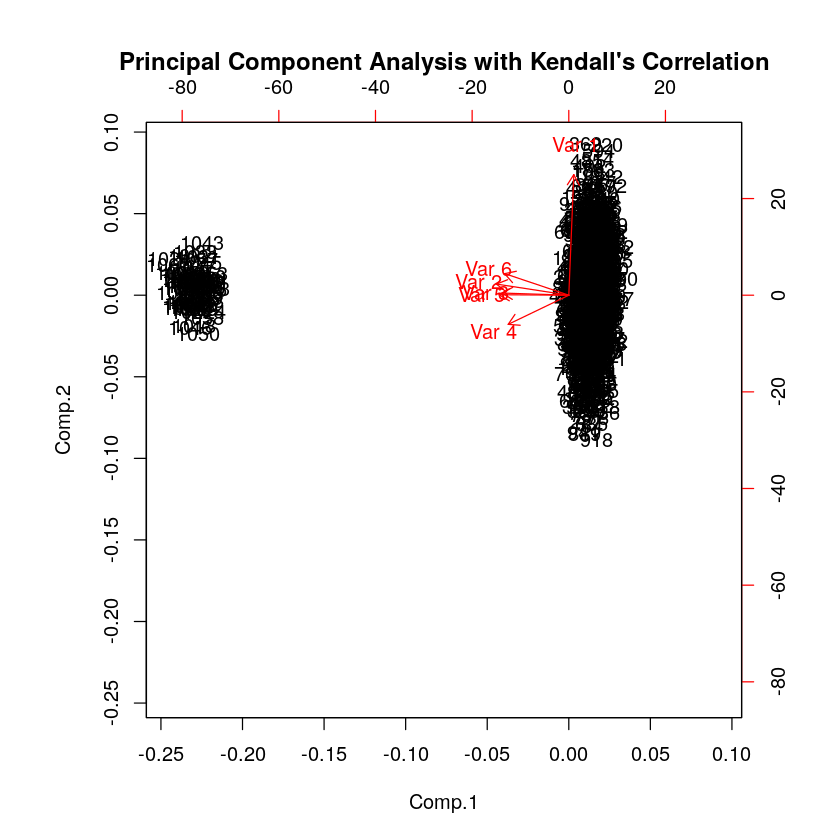}%
    \label{fig:plot5}}
  \subfloat[PCA using the Spearman's correlation for correlation matrix]{\includegraphics[scale = 0.35]{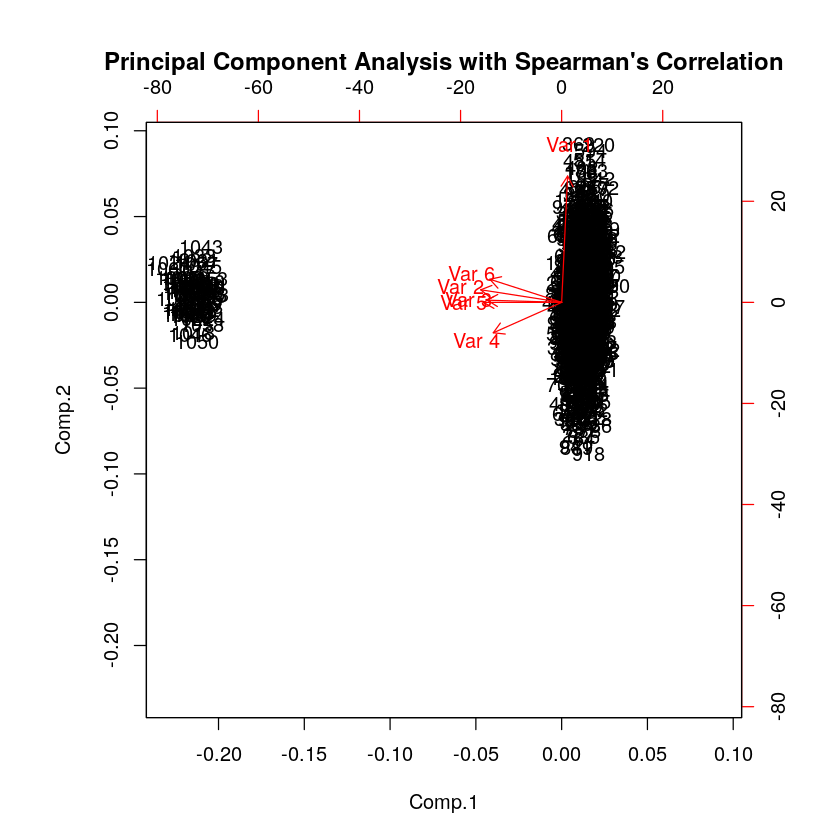}%
    \label{fig:plot6}}\\
    \subfloat[Classical PCA]{\includegraphics[scale = 0.35]{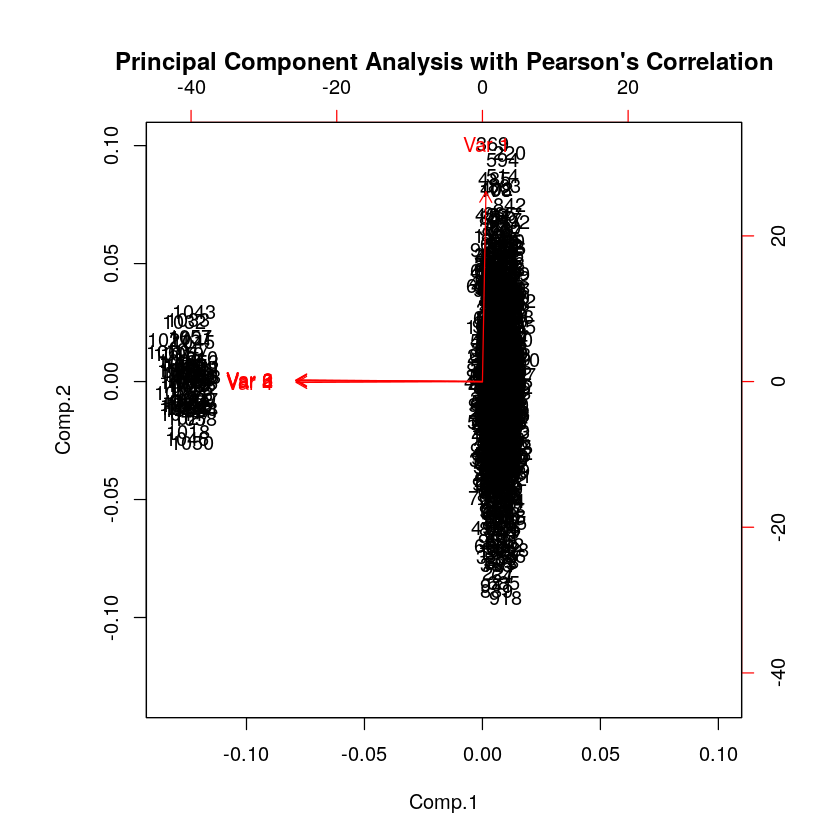}%
    \label{fig:plot7}}
    \caption{PCA in contaminated data with correlation matrix}\label{fig:correlation}
\end{figure}

\begin{figure}[ht!]
\centering
    \subfloat[PCA with robust covariance matrix]{\includegraphics[scale = 0.4]{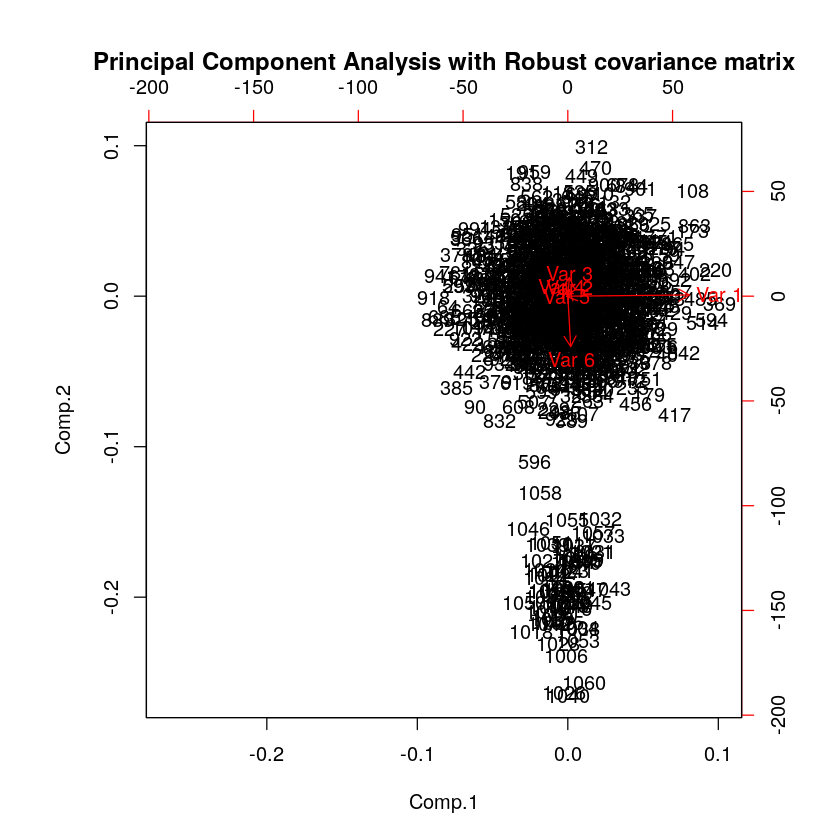}%
    \label{fig:plot8}}
  \subfloat[PCA with regular covariacne matrix]{\includegraphics[scale = 0.4]{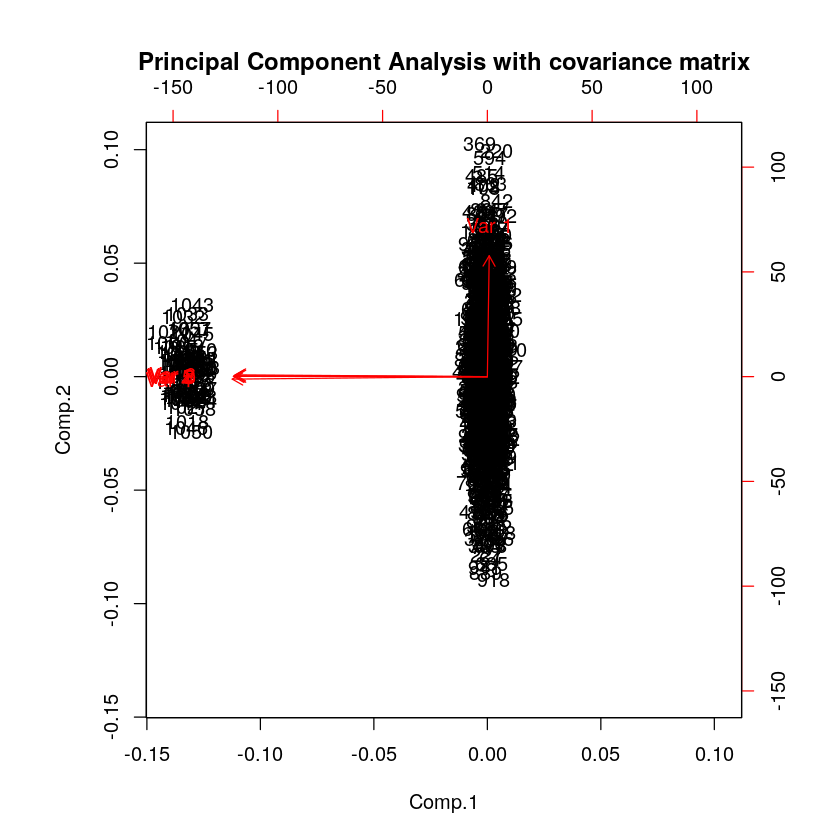}%
    \label{fig:plot9}}
    \caption{PCA in contaminated data with covariance matrix}\label{fig:covariance}
\end{figure}

Consider a data set of 6 normal distributed random variables with mean 
$0, 1, \dots, 5$ and covariance matrix given by diag$(5,1,1,1,1,1)$ with 
1000 observations each contaminated with 60 observations of a normal 
distribution with mean 20 and variance 5. In Figure~\ref{fig:correlation} 
are the graphics for the performance of the Principal Component Analysis 
with the correlation of Pearson (the regular correlation)(\ref{fig:plot7}),
the Kendall's $\tau$ (\ref{fig:plot5} and the Spearman's correlation 
(\ref{fig:plot6})\\[5pt]
The difference of performance between the last two no-parametric correlation 
coefficients between the variables is noticeable, in terms of how much the 
directions of the principal components respond to 
the presence of outliers in the data, compared to the regular PCA. It is also true
that the estimators tend to go a little to the contaminated data, but still
preserve a lot of the structure of the clean data.\\[5pt]
In Figure~\ref{fig:covariance} we can see the performance of the Principal 
component analysis using the robust estimation of the covariance matrix 
proposed in the previous section and the usual covariance matrix. Again, 
the usual is sensitive to outliers, while the robust matrix can perform 
better in the core of the data, given more precise directions for the principal 
components.

\subsection{Projection pursuit approach}
\begin{figure}[ht!]
\centering
    \subfloat[Classical PCA]{\includegraphics[scale = 0.25]{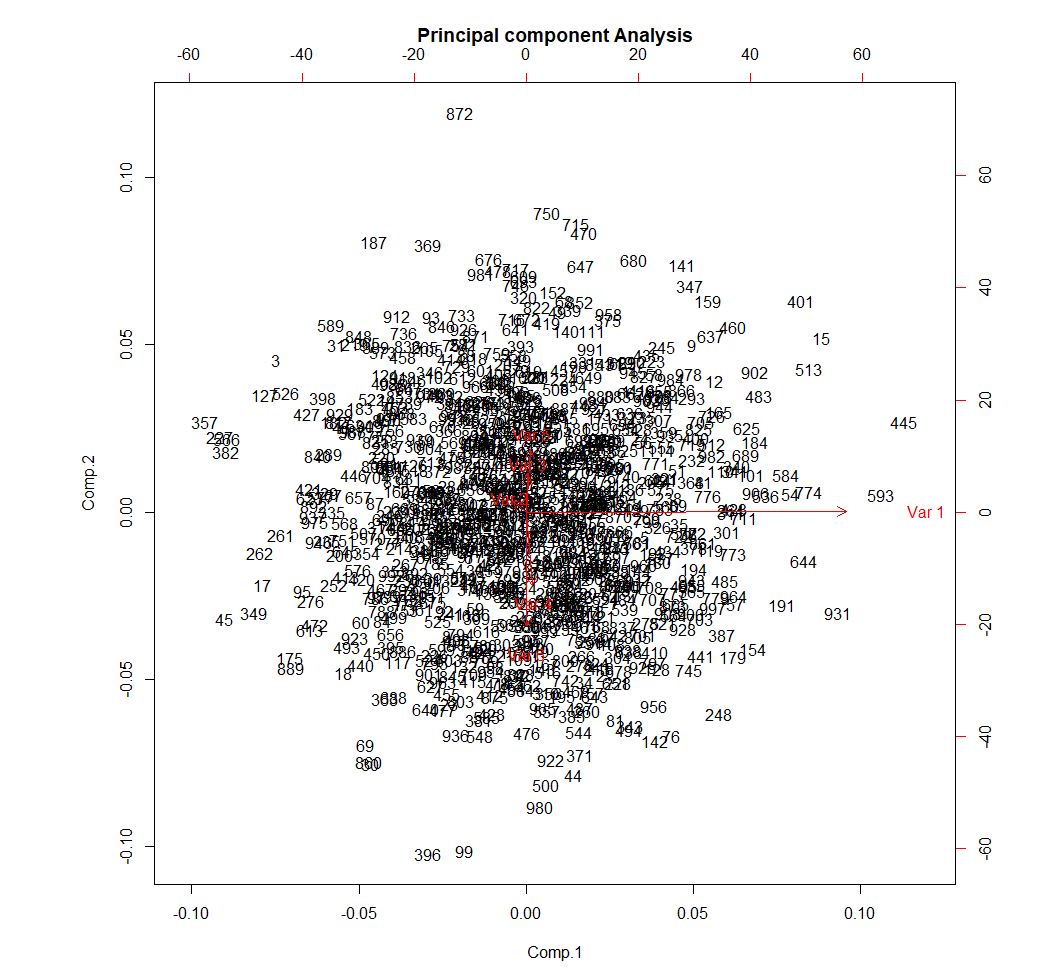}%
    \label{fig:plot2}}
  \subfloat[Projection Pursuit PCA]{\includegraphics[scale = 0.25]{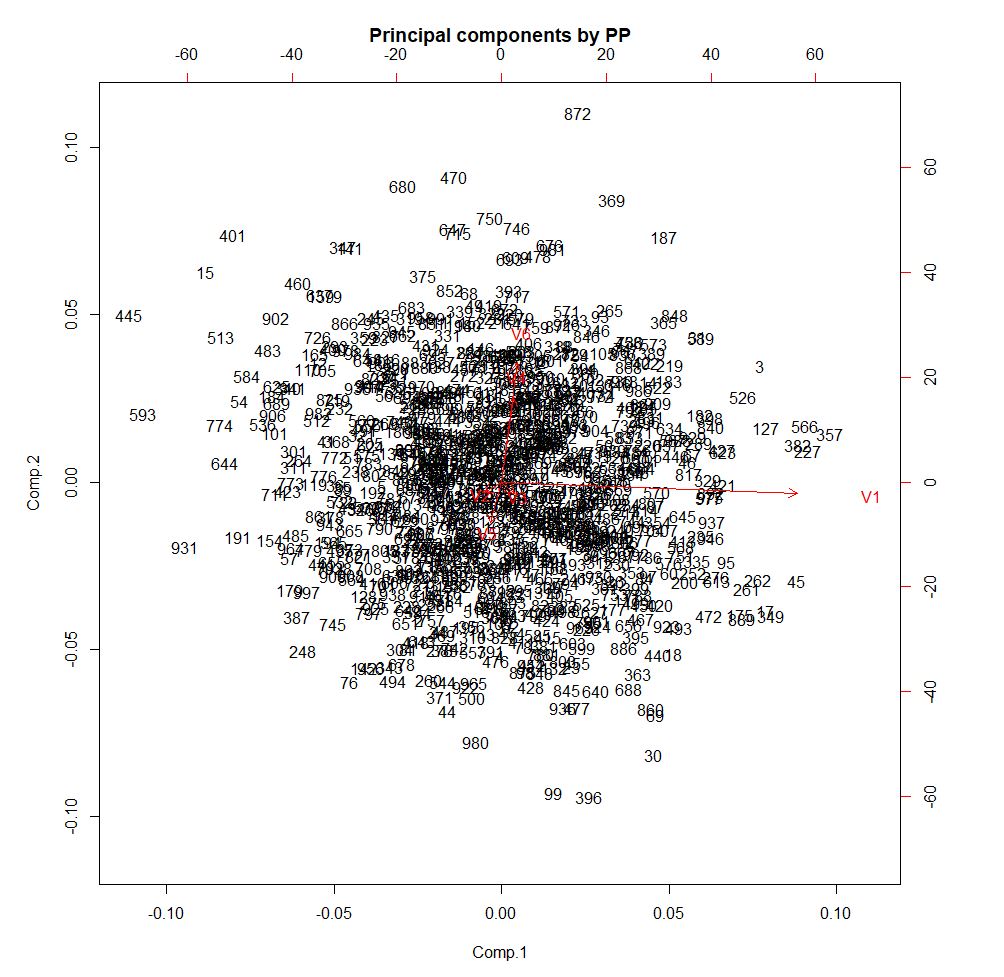}%
    \label{fig:plot1}}
    \caption{PCA in clean data}\label{fig:clean}
\end{figure}

Consider, again, a data set of 6 normal distributed random variables with mean 
$0, 1, \dots, 6$ and covariance matrix given by diag$(5,1,1,1,1,1)$ with 
1000 observations each. We apply the classical Principal Component Analysis 
to the data and obtain the results shown in the Figure~\ref{fig:plot2}, 
and also apply the algorithm for projection pursuit PCA, obtaining the
results shown in the Figure~\ref{fig:plot1}. As one can notice, the results are pretty similar. 

\begin{figure}[ht!]
\centering
    \subfloat[Classical PCA]{\includegraphics[scale = 0.25]{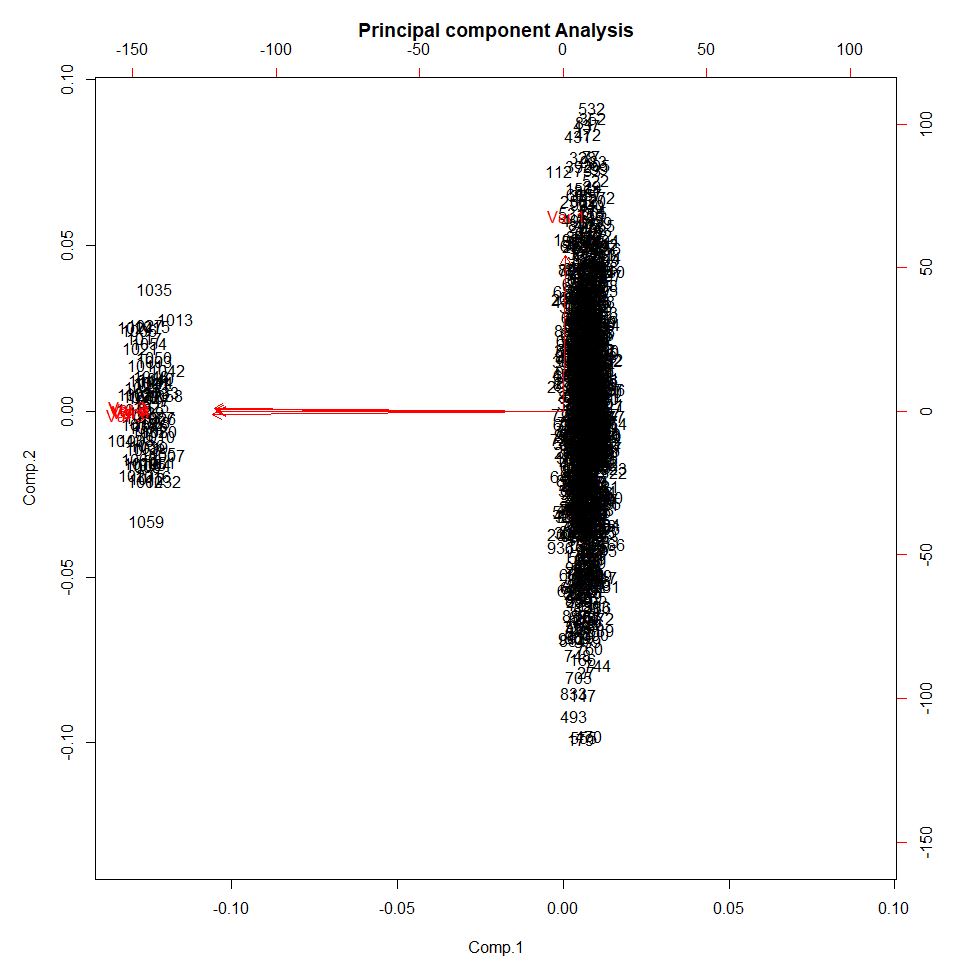}%
    \label{fig:plot4}}
  \subfloat[Projection Pursuit PCA]{\includegraphics[scale = 0.25]{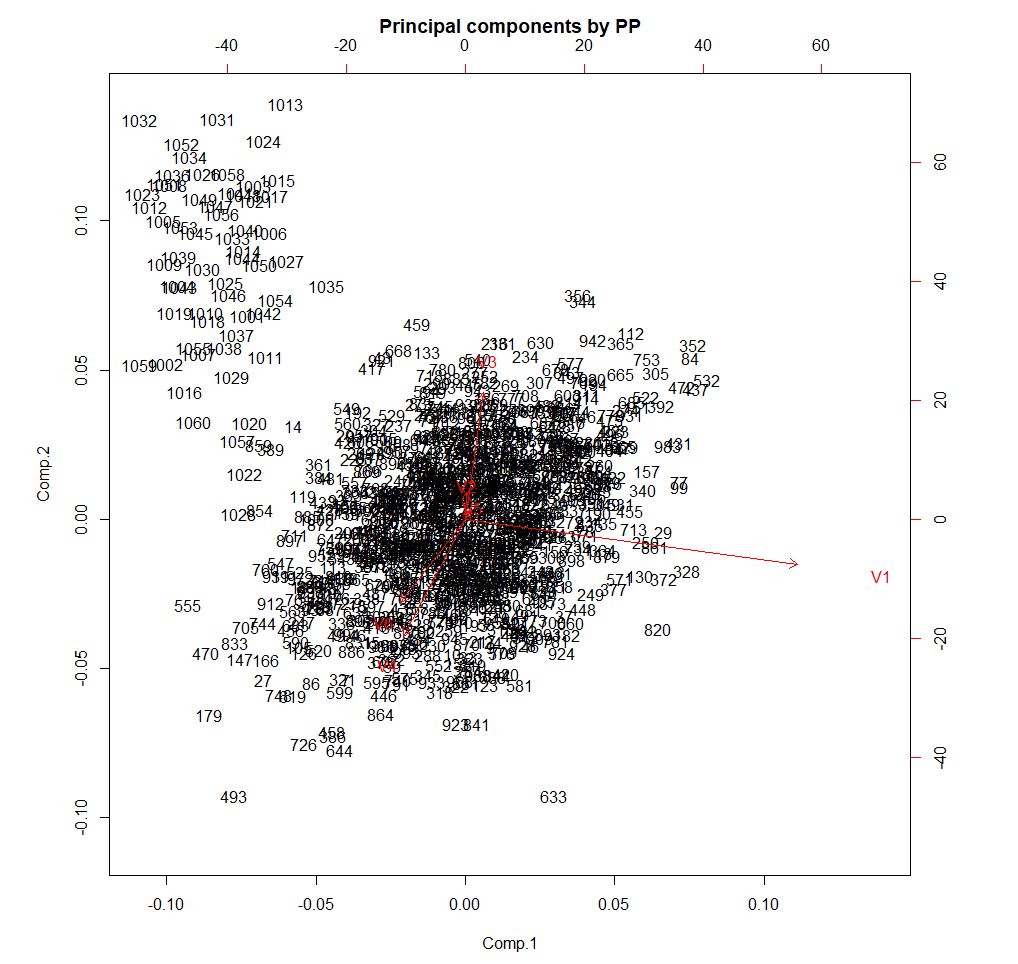}%
    \label{fig:plot3}}
    \caption{PCA in contaminated data}\label{fig:contaminated}
\end{figure}

Introducing atypical observations to the data, coming from a normal
distribution with mean 20 and variance 1, we analyze the performance of both methods in the
Figure~\ref{fig:contaminated}. Is easy to verify that classical PCA is very sensitive to analysis,
and the projection pursuit approach keeps the directions from the core of the data, surpassing
the outliers.

\section{Conclusions}
We presented three different ways to perform PCA that are in some
way robust to outliers presence in the data considered. One of them makes use of
robust estimators of theprevipus steps for the anlaysis itself, and the other
ones make the whole proces more robust to the presence of
atypical data.\\[5pt]
This is useful for the analysis of real problems, in which data of unknown distribution 
and possible contaminated values by human or numeric errors or wrong sampling can appear,
because it preserves a great part of the original structure of the data and gives points 
to analyze it from its core, without being drown to the outer layers by atypical values. 
{\small
\bibliographystyle{authordate1}

\bibliography{bibliografia}}

\end{document}